\documentclass[12pt]{article}

\usepackage{pstricks,pst-node}
\usepackage{psfig}
\usepackage{epsfig}
\usepackage{amsmath,amsthm}
\usepackage{amssymb}
\usepackage{color}
\usepackage{xspace}
\usepackage[colorlinks=true,
linkcolor=green,
filecolor=brown,
citecolor=green]{hyperref}

\def\v{\vert}

\def\a{\ensuremath{\mathcal A}\xspace}
\def\b{\ensuremath{\mathcal B}\xspace}
\def\c{\ensuremath{\mathcal C}\xspace}

\def\k{j}
\def\j{k}
\def\lra{left to right maxima\xspace}
\def\lrm{left to right maximum\xspace}

\def\mbf#1{\mathchoice{\hbox{\boldmath $\displaystyle #1$}}
        {\hbox{\boldmath $\textstyle #1$}}
        {\hbox{\boldmath $\scriptstyle #1$}}
        {\hbox{\boldmath $\scriptscriptstyle #1$}}} 

\newcommand{\StirlingPartition}[2]{\genfrac{ \{ }{ \} }{0pt}{}{#1}{#2}}

\newcommand{\seqnum}[1]{\href{http://oeis.org/#1}{\underline{#1}}}

\begin{document}
\newtheorem{theorem}{Theorem}
\newtheorem{defn}[theorem]{Definition}
\newtheorem{lemma}[theorem]{Lemma}
\newtheorem*{main}{Theorem}
\newtheorem{cor}[theorem]{Corollary}
\begin{center}
{\Large
A permutation pattern that illustrates the strong law of small \\[2mm] numbers    
}

\vspace*{5mm}

DAVID CALLAN  \\
\noindent {\small Dept. of Statistics, 
University of Wisconsin-Madison,  Madison, WI \ 53706}  \\
{\bf callan@stat.wisc.edu} 
\end{center}

\begin{abstract}
We obtain an explicit formula for the number of permutations of $[n]$ 
that avoid the barred pattern $\bar{1}43\bar{5}2$. A curious feature 
of its counting sequence, $1, 1, 2, 5, 14, 43,$ $145, 538, 
2194,\ldots$, is that the displayed terms agree with A122993 in the 
On-Line Encyclopedia of Integer Sequences, but the two sequences 
diverge thereafter.  

\end{abstract}


\section{Introduction}
\vspace*{-4mm}
A permutation $\pi$ avoids the barred pattern $\bar{1}43\bar{5}2$
if each instance of a not-necessarily-consecutive 432 pattern in $\pi$ 
is part of a 14352 pattern in $\pi$, and similarly for other barred patterns. 
This paper is one of a series of notes counting permutations avoiding 
a 5 letter pattern with 2 bars that do not yield to Lara Pudwell's 
method of Enumeration Schemes \cite{schemes}. The question of whether there may be an automated method to fill in these and other gaps in Pudwell's enumeration remains open. Here we treat the 
pattern $\bar{1}43\bar{5}2$.
A curious feature 
of the counting sequence is that it agrees through the $n=8$ term  with 
sequence \seqnum{A122993} in the 
On-Line Encyclopedia of Integer Sequences \cite{oeis}, an instance of 
the Strong Law of Small Numbers \cite{strong1,strong2}.

Our method is to identify the structure of a 
$\bar{1}43\bar{5}2$-avoider. This permits a direct count 
as a 5-summation formula according to five statistics of the permutation, 
four of which are the first entry $a$, the immediate predecessor of 1 denoted 
$b$, the position of 1 denoted $\k$, and the number of left to right 
maxima that occur after 1 denoted $\j$. One of these sums can be evaluated, leading to a faster 
formula.

\vspace*{-4mm}

\section[bar\{1\}43bar\{5\}2-Avoiders]{$\mbf{\bar{1}43\bar{5}2}$-Avoiders}
\vspace*{-3mm}
A ``typical'' $\bar{1}43\bar{5}2$-avoider is illustrated in Figure 1 in matrix 
form.  It has first entry $a=5$,  
1 is in position $\k=4$, the immediate predecessor of 1 is $b=16$, 
and there are $\j=5$ left to right 
maxima that occur after 1. Here $\k\ge 3$ so that $1,a,b$ are all 
distinct. The special cases $\k=1$ or 2 are treated later. There is a 
vertical blue line through the bullet representing the entry 1, and 
yellow vertical lines through the \lra that occur after 1. These $\j$ 
yellow lines divide the the part of the matrix to the right of 
the blue line into $\j+1$ vertical strips (in white, one of which is 
vacuous in Figure 1). Furthermore, horizontal lines through 1, $a$ 
and $b$ determine three horizontal strips indexed by 
$\a=[2,a-1],\ \b=[a+1,b-1],\ 
\c=[b+1,n]$.
There are $\k-3$ bullets to the left of the 
blue line in strip $\b$ and none in $\a$ or $\c$. Hence, to the right of the blue line there are 
$A:= a-2$ bullets in strip $\a$, $B:= \v\,\b\,\v-(\k-3) =
b-a-\k+2$ bullets in strip $\b$, and  $C:=
n-b$ bullets in strip $\c$.
\begin{figure}
\vspace{-.6in}
\begin{center}
\includegraphics[angle=0, scale = 1]{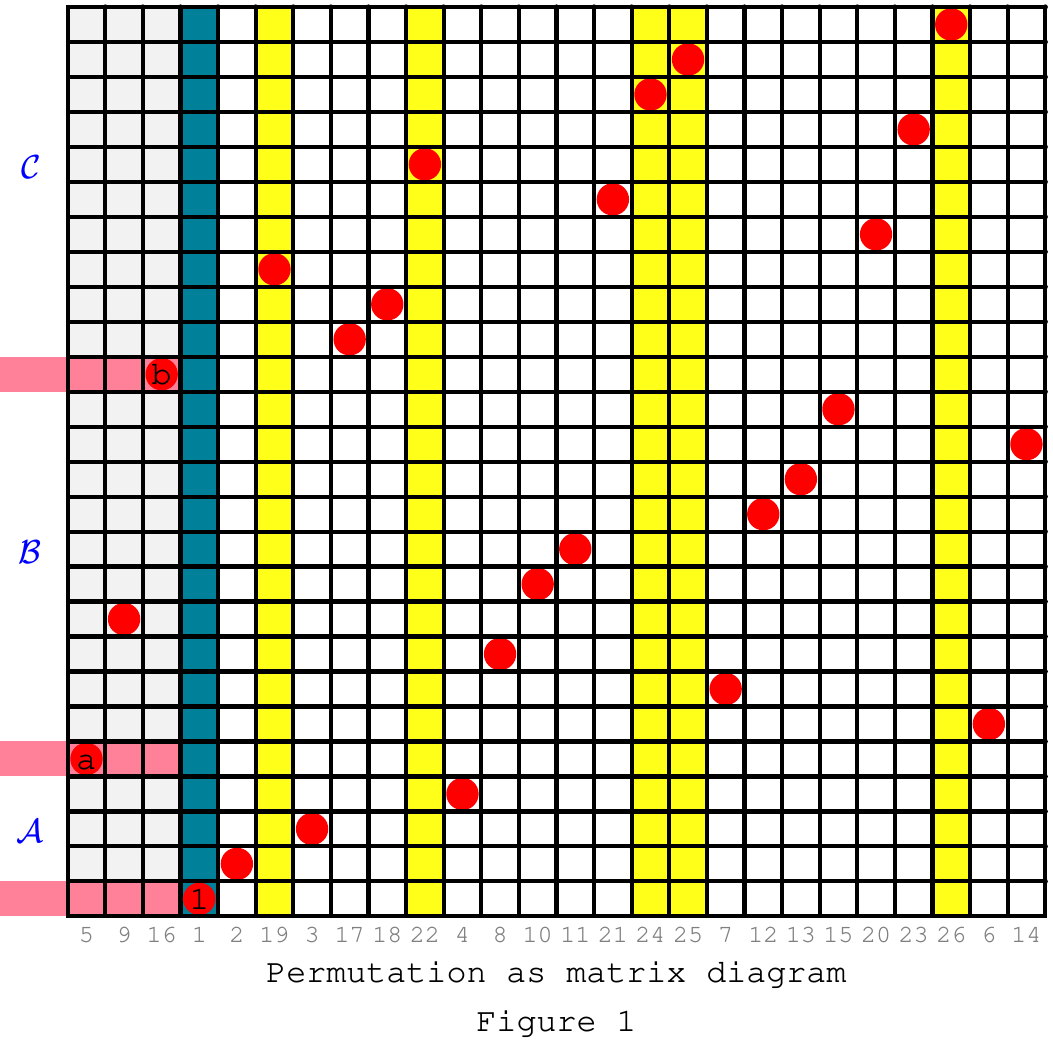}
\end{center}
\end{figure}
The following properties of a $\bar{1}43\bar{5}2$-avoider are evident 
in the illustration and easily proved from the definition.

\vspace*{-3mm}
\begin{itemize}
    \item  The entries to the left of 1 are increasing, 
    else together with 1, there is a 432 pattern with no 
    available 1. Equivalently, the bullets in the gray vertical 
    strip on the left are rising.
    \vspace*{-2mm}
    \item  The entries $2,3,\ldots,a-1$ occur in that order, else together 
    with $b$, there is a 432 pattern with no available 1. Equivalently, 
    the bullets in horizontal strip $\a$ are rising.
    \vspace*{-2mm}
    \item Entries in the interval $(a,b)$ lie either to the left of 1 or to the 
    right of $a-1$, else together  with $b$, 
    there is a 432 pattern with no available 1. Equivalently, 
    all bullets in horizontal strip $\b$ to the right of 1 are also to 
    the right of $a-1$.
    \vspace*{-2mm}
    \item  Every descent initiator after 1 is a \lrm, else together with a 
    \lrm to its left (there is one), we have a 432 pattern with no available 5.
     Equivalently, the bullets in each vertical white strip $\a$ are rising.
\end{itemize}
\vspace*{-2mm}
Conversely, when the position $\k$ of 1 is $\ge 3$, one can check 
that a permutation with these 
properties is  $\bar{1}43\bar{5}2$-avoiding. 

To count permutations with these four properties, let $i\in [1,\j+1]$ denote the left to 
right position of the first white strip containing an entry in 
$(a,b)$, that is, containing a bullet in horizontal strip $\b$ (when 
there is one).

The subpermutation of entries in $\c=[b+1,n]$, when split at its 
\lra, forms a partition in 
a canonical form: in each block, the largest entry occurs first and the 
rest of the block is increasing, and the blocks are ordered by increasing 
first entry. This yields $\StirlingPartition{C}{\j}$ choices to 
determine the relative positions of entries in $\c$. 

Next, choose $\k-3$ elements 
from $[a+1,b-1]$ to precede 1---$\binom{b-a-1}{\k-3}$ choices. The entries 
in $\b$ following 1 
must be distributed into boxes (white strips) labeled 
$i,i+1,\ldots,\j+1$ in such a way that box $i$ is 
nonempty---$(\j-i+2)^{B}-(\j-i+1)^{B}$ choices when $B>0$. The bullets for 
entries in $\a$ must be distributed into boxes 
$1,2,\ldots,i$---$\binom{A+i-1}{i-1}$ choices when $B>0$. In case 
$B=0$, we merely distribute the bullets for entries in $\a$ into $\j+1$ 
boxes---$\binom{A+\j}{\j}$ choices. 

Recalling that $A=a-2,\ B=b-a-\k+2,\ C=n-b$, the contribution of the case 
$\k\ge 3$ to the desired count is now seen to be
\begin{multline}\label{base}
\sum_{a=2}^{n-1}\sum_{b=a+1}^{n}\sum_{\k=3}^{b-a+1} 
\sum_{\j=1}^{n-b}\sum_{i=1}^{\j+1}
  \StirlingPartition{n-b}{\j} \binom{b-a-1}{\k-3}
 \big( (\j-i+2)^{b-a-\k+2}-(\j-i+1)^{b-a-\k+2} \big) \times \\
  \binom{a+i-3}{i-1} +
\sum_{a=2}^{n-1}\sum_{b=a+1}^{n} 
  \StirlingPartition{n-b}{\j} 
 \binom{a+i-3}{i-1}
\end{multline}

When $\k=1$, the map ``delete first entry'' is a bijection to 
$43\bar{5}2$-avoiding permutations of size $n-1$, counted by the 
Bell number $B_{n-1}$ \cite{eigensequence}. When $\k=2$, we have $a=b$ in Figure 1, and the 
count reduces to $\sum_{a=2}^{n}\sum_{\j=0}^{n-a}\binom{\j+a-2}{a-2}\StirlingPartition{n-a}{\j}$ where 
$\StirlingPartition{0}{0}:=1$.

The sum over $\k$ in (\ref{base}) can be evaluated using the binomial theorem, and putting it all together we have, after minor simplifications, the following result.
\begin{main}
For $n\ge 2$, the number of permutations of $[n]$ avoiding the barred pattern $\bar{1}43\bar{5}2$ is
\begin{multline*}
B_{n-1} + 1 + 2^{n-2}-n\ + \\
\sum_{a=0}^{n-3}\sum_{b=0}^{a-1}\sum_{\j=0}^{a-b} \left(\sum_{i=0}^{\j}
\binom{n - 4 - a + \j - i}{\j - i}(i+2)^b - \binom{n - 3 - a + \j}{\j}\right) \StirlingPartition{a - b}{\j}\: + \\
\sum_{a=0}^{n-2}\sum_{\j=0}^{n-2-a}\binom{\j + a + 1}{\j + 1}\StirlingPartition{n-2-a}{\j}.
\end{multline*}
\end{main}
The first few terms of the counting sequence, starting at $n=1$, are 1,\:2,\:5,\:14,\:43,\:145,
\newline 
\:538,\:2194,\:9790,\:47491,\:248706.


\begin{thebibliography}{99}


\bibitem{invert} David Callan, 
The number of $\bar{2}413\bar{5}$-avoiding permutations, preprint,
\htmladdnormallink{arXiv:1110.6884}{http://front.math.ucdavis.edu/1110.6884}, 4pp.

\bibitem{pudwell} David Callan, 
The number of $\bar{3}\bar{1}542$-avoiding permutations, preprint,
\htmladdnormallink{arXiv:1111.3088}{http://front.math.ucdavis.edu/1111.3088}, 5pp.
    
\bibitem{schemes}  Lara Pudwell, Enumeration Schemes for Permutations Avoiding 
Barred Patterns, 
\htmladdnormallink{\emph{Electronic J. Combinatorics}}{http://www.combinatorics.org/}
 \textbf{17} (1) (2010), R29, 27 pp.

\bibitem{oeis}
The On-Line Encyclopedia of Integer Sequences, published electronically at 
\htmladdnormallink{http://oeis.org}{http://oeis.org}, 2010.


\bibitem{strong1} Guy, Richard K.,
The Strong Law of Small Numbers, \emph{American Mathematical Monthly} 
\textbf{95} Issue 8 (1988), 697-Ð712.

\bibitem{strong2} Guy, Richard K.,
The Second Strong Law of Small Numbers, \emph{Mathematics Magazine} 
\textbf{63} (1990), 3--20.

\bibitem{eigensequence} David Callan, A combinatorial interpretation of the eigensequence for composition,
\htmladdnormallink{\emph{J. Integer Sequences}}{http://www.cs.uwaterloo.ca/journals/JIS/} \textbf{9} (2006) Article 06.1.4, 12 pp. 
 

\end{thebibliography}
\end{document}